\documentclass[12pt,a4paper]{article}

\usepackage{epstopdf}
\usepackage{array}
\usepackage{dcolumn,anysize}
\usepackage{graphicx}
\usepackage{float}
\usepackage{enumerate,color}
\usepackage{graphics}
\usepackage{amsmath,amsfonts,amsthm}
\usepackage[cp850]{inputenc}

\newcommand{\Eqref}[1]{(\ref{#1})}

\title{Analytic formulas for the evaluation of the Pearcey integral}
\author{\\Jos\'{e} L. L\'{o}pez$^1$ and Pedro J. Pagola$^1$\\\\
\small{$^1$
\textsf{\textit{Dpto. de Ingenier\'{\i}a Matem\'{a}tica e Inform\'{a}tica, Universidad P\'{u}blica de Navarra}}}\\
\small{ \textsf{\textit{ e-mail: jl.lopez@unavarra.es, pedro.pagola@unavarra.es}}}}
\date{}

\begin{document}
\normalsize \maketitle

\begin{abstract}
We can find in the literature several convergent and/or asymptotic expansions of the Pearcey integral $P(x,y)$ in different regions of the complex variables $x$ and $y$, but they do not cover the whole complex $x$ and $y$ planes. The purpose of this paper is to complete this analysis giving new convergent and/or asymptotic expansions that, together with the known ones, let the evaluation of the Pearcey integral in a large region of the complex $x$ and $y$ planes. The accuracy of the approximations derived in this paper is illustrated with some numerical experiments. Moreover, the expansions derived here are simpler compared with other known expansions, as they are derived from a simple manipulation of the integral definition of $P(x,y)$.
 \\\\
\noindent \textsf{2010 AMS \textit{Mathematics Subject
Classification:} 33E20; 41A60.
} \\\\
\noindent  \textsf{Keywords \& Phrases:} Pearcey integral. Convergent and asymptotic expansions. Watson lemma.
\\\\
\end{abstract}

\section{Introduction}

The Pearcey integral is the second canonical diffraction integrals \cite[Sec. 36.14]{nist}. Apart from their mathematical importance in the uniform asymptotic approximation of oscillatory integrals \cite{olde}, the canonical diffraction integrals have physical applications in the description of surface gravity waves \cite{kelvin}, \cite{ursell}, bifurcation sets, optics, quantum mechanics and acoustics (see \cite[Sec. 36.14]{nist} and references there in). The Cusp catastrophe or Pearcey integral \cite[p.777, eq. 36.2.14]{nist} is the integral:
\begin{equation}\label{pea}
\bar P(x,y):=\int_{-\infty}^\infty e^{i(t^4 +xt^2+yt)}dt.
\end{equation}
This integral was first evaluated numerically by using quadrature formulas in \cite{pearcey} in the context of the investigation of the electromagnetic field near a cusp. In \cite[Chap. 36]{nist} we can find many properties such as symmetries, ilustrarive pictures, bifurcation sets, scaling relations, zeros, convergent series expansions, differential equations and leading-order asymptotic approximations among others. But we cannot find many details about analytic approximation formulas (asymptotic expansions in particular) or numerical evaluation techniques.

The integral \Eqref{pea} exists only for $0\le\arg{x}\le \pi$ and real $y$. As it is indicated in \cite{paris}, after a rotation of the integration path through an angle of $\pi/8$ that removes the rapidly oscillatory term $e^{it^4}$, the Pearcey integral may be written in the form $\bar P(x,y)=2e^{i\pi/8}P(xe^{-i\pi/4},ye^{i\pi/8})$, with
\begin{equation}\label{pearcy}
P(x,y):=\int_0^\infty e^{-t^4 -xt^2}\cos(yt)dt=\frac{1}{2}\int_{-\infty}^{\infty} e^{-t^4+iyt}e^{-xt^2}dt.
\end{equation}
This integral is absolutely convergent for all complex values of $x$ and $y$ and represents the analytic continuation of the Pearcey integral $\bar P(x,y)$ to all complex values of $x$ and $y$ \cite{paris}. Therefore, it is more convenient to work with the representation \Eqref{pearcy} of the Pearcey integral.

We summarize below the more elementary analytic expansions (convergent or asymptotic) of $P(x,y)$ that we can find in the literature and give references to other more elaborated expansions.

A convergent series expansion of $P(x,y)$ may be found in \cite[p. 787, eqs. 36.8.1 and 36.8.2]{nist}:
\begin{equation}\label{nistexpan}
P(x,y)=\frac{1}{4}\sum_{n=0}^\infty (-1)^n\Gamma\left(\frac{2n+1}{4}\right)a_{2n}(x,y),
\end{equation}
with $a_0(x,y)=1$, $a_1(x,y)=y$ and, for $n=2,3,4,...$,
$$
a_n(x,y)=\frac{1}{n}[y\,a_{n-1}(x,y)+2x\,a_{n-2}(x,y)].
$$
Another convergent expansion of $P(x,y)$ may be obtained after an expansion of the cosine term in \Eqref{pearcy} and interchange of sum and integral \cite{pariskaminski}:
\begin{equation}\label{confluexpan}
P(x,y)=\sum_{n=0}^\infty \frac{(-y^2)^n}{(2n)!}P_n(x), \hskip 2cm (x,y)\in \mathbb{C}^2,
\end{equation}
with
\begin{equation}\label{eme}
P_n(x):=\left\lbrace
\begin{array}{cl}
\frac{1}{2^{n+3/2}}\Gamma\left(n+\frac{1}{2}\right)U\left(\frac{n}{2}+\frac{1}{4},\frac{1}{2};\frac{x^2}{4}\right) & \hskip 6cm {\rm if}\hskip 2mm \Re x\ge 0, \\\\
\frac{1}{4}\Gamma\left(\frac{n}{2}+\frac{1}{4}\right)M\left(\frac{n}{2}+\frac{1}{4},\frac{1}{2};\frac{x^2}{4}\right)-&
\frac{x}{4}\Gamma\left(\frac{n}{2}+\frac{3}{4}\right)M\left(\frac{n}{2}+\frac{3}{4},\frac{3}{2};\frac{x^2}{4}\right) \hskip 5mm {\rm if}\hskip 2mm \Re x< 0,
\end{array}
\right.
\end{equation}
where $M(a,b;z)$ and $U(a,b;z)$ are confluent hypergeometric functions \cite[Chap. 13]{conflus}

Apart from convergent expansions, we can also find in the literature several asymptotic expansions of $P(x,y)$.
In \cite{kaminski} we can find an asymptotic expansion of the Pearcey integral when $(x,y)$ are near the caustic $8x^3-27y^2=0$ that remains valid as $\vert x\vert\to\infty$. The expansion is given in terms of Airy functions and its derivatives and the coefficients are computed recursively. We refer the reader to \cite{kaminski} for further details.

An exhaustive  asymptotic analysis of this integral can be found in \cite{paris}. In particular, a complete asymptotic expansion is given in \cite{paris} by using asymptotic techniques for integrals applied to the integral \Eqref{pearcy}. The asymptotic analysis of this integral for large $\vert x\vert$ is divided in two regions: $\vert\arg{x}\vert<\frac{\pi}{2}$ and $\vert\arg{x}\vert>\frac{\pi}{2}$. In the first region we find that \cite{paris}
\begin{equation}\label{parisi}
P(x,y)\sim\frac{1}{2\sqrt{x}}e^{-y^2/(4x)}\sum_{n=0}^\infty\frac{(-1)^n\Gamma(2n+1/2)}{n!x^{2n}}M\left(-2n;1/2;\frac{y^2}{4x}\right).
\end{equation}
The asymptotic expansion in the second region is a little bit more cumbersome, we refer to \cite{paris} for details. Another convergent and asymptotic expansion of $P(x,y)$ for large $\vert x\vert$, derived from a differential equation satisfied by $P(x,y)$, may be found in \cite{pedrox}.

A complete asymptotic expansion for large $\vert y\vert$ is derived in \cite{pedro}, where the asymptotic analysis is divided in three different regions: $\vert\arg y\vert\le\frac{\pi}{8}$, $-\frac{\pi}{2}\le\arg y<-\frac{\pi}{8}$ and $\frac{\pi}{8}<\arg y\le\frac{\pi}{2}$. We give details here for the last region and refer the reader to \cite{pedro} for further details. When $\frac{\pi}{8}<\arg y\le\frac{\pi}{2}$ we have that:
\begin{equation}\label{cucu}
P(x,y)\sim \sqrt{\frac{\pi}{3}}\frac{1}{2^{5/6}y^{1/3}}e^{3\cdot 4^{-4/3}e^{-2i\pi/3}y^{4/3}-4^{-2/3}xy^{2/3}e^{-i\pi/3}+x^2/6+i\pi/6}\sum_{n=0}^\infty\frac{A_n(x)}{y^{2n/3}},
\end{equation}
with
$$
A_n(x):=\displaystyle{\sum_{m=\lfloor(n+1)/2\rfloor}^n\sum_{k=0}^{2m-n}b_{n,m,k}(x)c_{n+2m-k}(x)},
$$
where
$$b_{n,m,k}(x):=\frac{x^k\,2^{4(2m-n-k)/3}\,(-1)^{m}}{k!(2m-n-k)!(n-m)!}$$
and
$$
c_{n+2}(x)=\displaystyle{\frac{x}{3\cdot 2^{1/3}}}c_{n+1}(x)+\displaystyle{\frac{n+1}{3\cdot 2^{2/3}}}c_n(x),\hskip 2cm n=0,1,2,...
$$
with $c_0(x)=1$ and $c_1(x)=\displaystyle{\frac{x}{3\cdot 2^{1/3}}}$.

In \cite{pariskaminski} we can find the hyperasymptotic evaluation of the Pearcey integral for real values of $x$ and $y$ using Hadamard expansions. The Pearcey integral is written in terms of an infinite series whose terms are Hadamard series. The terms of these Hadamard series are incomplete gamma functions. The analytic expression is sophisticated and we refer the reader to \cite{pariskaminski} for details. In this paper we exploit the following simple idea: the integrand in \eqref{pearcy} is an exponential. Then, we manipulate the function in the exponent in order to factorize the exponential and replace one of the factors by its Taylor expansion at the origin. We consider three different possibilities.

In the following section we derive a complete convergent and asymptotic expansion of $P(x,y)$ for small $\vert x\vert$. In Section 3 we derive a complete asymptotic expansion of $P(x,y)$ for large $\vert x\vert$. In Section 4 we derive a complete asymptotic expansion of $P(x,y)$ for large $\vert x\vert $ and $\vert y\vert$.  Section 5 contains some numerical experiments and a few remarks. Through the paper we use the principal argument $\arg z\in(-\pi,\pi]$ for any complex number $z$.

\section{A convergent and asymptotic expansion for small $x$}

Because $P(x,y)=P(x,-y)$, without loss of generality, we may restrict ourselves to $\Re y\ge 0$. A convergent expansion of $P(x,y)$ that is also and asymptotic expansion when $\vert x\vert\to 0$ can be obtained by considering the Taylor expansion of the exponential at the origin:
\begin{equation}\label{expanuno1}
e^{-xt^2}=\sum_{k=0}^{n-1}\frac{(-x)^kt^{2k}}{k!}+r_n\left(xt^2\right),
\end{equation}
where $r_n(t)$ is the Taylor remainder. Replacing this expansion in the integrand of the right hand side of formula \eqref{pearcy} and interchanging sum and integral we find:
\begin{equation}\label{pe}
P(x,y)=\displaystyle{\sum_{k=0}^{n-1} \frac{(-x)^k}{k!}P_{k}(y)+R_n(x,y),}
\end{equation}
with
$$
P_k(y):=\frac{1}{2}\int_{-\infty }^{\infty } t^{2k}e^{-t^4+iy t}dt,\hskip 2cm R_n(x,y):=\frac{1}{2}\int_{-\infty }^{\infty } r_n\left(xt^2\right)e^{-t^4+iy t}dt.
$$
After some straightforward manipulations, we obtain an explicit formula for the coefficients $P_k(y)$ in \eqref{pe} $\forall k=0,1,2,\ldots$:
$$
P_k(y):=\frac{1}{4}\Gamma\left(\frac{1+2k}{4}\right){}_1F_3\left(\frac{1+2k}{4};\frac{1}{4},\frac{2}{4},\frac{3}{4};\left(\frac{y}{4}\right)^4\right)-\frac{y^2}{8}\Gamma\left(\frac{3+2k}{4}\right){}_1F_3\left(\frac{3+2k}{4};\frac{3}{4},\frac{5}{4},\frac{6}{4};\left(\frac{y}{4}\right)^4\right).
$$

Integrating by parts we derive the following five terms recurrence for the coefficients $P_k(y)$:
$$
P_{k+4}(y)=-\left(\frac{y}{4}\right)^2\, P_{k+1}(y)-\frac{(k+1)(2k+1)}{8}P_k(y)+(k+1)(k+7/4)P_{k+2}(y)\hskip 0.5cm, \hskip 0.5cm k=0,1,2,....
$$

with

$$
P_0(y)=\frac{1}{4}\Gamma\left(\frac{1}{4}\right){}_0F_2\left(-;\frac{1}{2},\frac{3}{4};\frac{y^4}{256}\right)-\frac{y^2}{8}\Gamma\left(\frac{3}{4}\right){}_0F_2\left(-;\frac{5}{4},\frac{3}{2};\frac{y^4}{256}\right),
$$
$$
P_1(y)=\frac{1}{4}\Gamma\left(\frac{3}{4}\right){}_0F_2\left(-;\frac{1}{2},\frac{1}{4};\frac{y^4}{256}\right)-\frac{y^2}{8}\Gamma\left(\frac{5}{4}\right){}_0F_2\left(-;\frac{3}{4},\frac{3}{2};\frac{y^4}{256}\right),
$$
$$
P_2(y)=\frac{1}{4}\Gamma\left(\frac{5}{4}\right){}_1F_3\left(\frac{5}{4};\frac{1}{4},\frac{1}{2},\frac{3}{4};\frac{y^4}{256}\right)-\frac{y^2}{8}\Gamma\left(\frac{7}{4}\right){}_1F_3\left(\frac{7}{4};\frac{3}{4},\frac{5}{4},\frac{3}{2};\frac{y^4}{256}\right),
$$
$$
P_3(y)=\frac{1}{4}\Gamma\left(\frac{7}{4}\right){}_1F_3\left(\frac{7}{4};\frac{1}{4},\frac{1}{2},\frac{3}{4};\frac{y^4}{256}\right)-\frac{y^2}{8}\Gamma\left(\frac{9}{4}\right){}_0F_2\left(\frac{9}{4};\frac{3}{4},\frac{5}{4},\frac{3}{2};\frac{y^4}{256}\right).
$$

Therefore, we have obtained that, for any complex $x$ and $y$,
\begin{equation}\label{expanuno}
P(x,y)=\displaystyle{\sum_{k=0}^{n-1} \frac{P_{k}(y)}{k!}(-x)^k+R_n(x,y).}
\end{equation}
This expansion is convergent for any value of $x$ and $y$ and is an asymptotic expansion of $P(x,y)$ as $x\to 0$.
Moreover, using the Lagrange formula for the remainder in \eqref{expanuno1}
$$
r_n\left(xt^2\right)=\frac{e^{-\xi\,\Re x}}{n!}(xt^2)^n\,\,\,\,\textnormal{with}\,\,\,0<\xi<t^2
$$
we have

$$
\vert r_n\left(xt^2\right)\vert\,<\,\left\lbrace
\begin{array}{ll}
\frac{t^{2n}}{n!}\vert x\vert^n & \hskip 5mm {\rm if}\hskip 2mm \Re x\ge 0, \\\\
\frac{t^{2n}e^{\vert \Re x\vert\,t^2}}{n!}\vert x\vert^n& \hskip 5mm {\rm if}\hskip 2mm \Re x< 0,
\end{array}
\right.
$$
In the first case ($\Re x\ge 0$), we immediately conclude that
$$
\vert R_n(x,y)\vert\le\frac{\vert x\vert^n}{n!}P_n\left(i\,\Im y\right).
$$
In second case ($\Re x< 0$), we have
$$
\vert R_n(x,y)\vert=\left\vert\frac{1}{2}\int_{-\infty }^{\infty } r_n\left(xt^2\right)e^{-t^4+iy t}dt\right\vert\leq \frac{1}{2}\frac{\vert x\vert^n}{n!}\int_{-\infty }^{\infty} t^{2n}e^{\vert \Re x\vert\,t^2}e^{-t^4-\Im y t}dt.
$$
Writing the dominant term $e^{-t^4}$ of integrand in the form $\displaystyle e^{-\frac{t^4}{2}}\,e^{-\frac{t^4}{2}}$ we can rewrite the above integral as
$$
\frac{1}{2}\frac{\vert x\vert^n}{n!}\int_{-\infty }^{\infty}g(t)\cdot e^{-\frac{t^4}{2}-\Im y t}t^{2n}dt\hskip 1cm\textnormal{with}\hskip 1cm g(t):=e^{-\frac{t^4}{2}+\vert \Re x\vert\,t^2}.
$$
The function $g(t)$ attains it absolute maximum at $t=\sqrt{\vert \Re x\vert}$ and then $\displaystyle g(t)\leq g(\sqrt{\vert \Re x\vert})=e^{\frac{\Re^2(x)}{2}}$. Therefore, we conclude that, for $\Re x<0$,
$$
\vert R_n(x,y)\vert\le\frac{\vert x\sqrt{2}\vert^n}{n!}\,\sqrt[4]{2}\,e^{\frac{ \Re^2(x)}{2}}\,P_n\left(i\,\sqrt[4]{2}\,\Im y\right).
$$

\section{An asymptotic expansion for large $\vert x\vert$}

Using the Cauchy's residue theorem we find that we can replace the integration path $(-\infty,\infty)$ in \Eqref{pearcy} by any straight $\Gamma:=\lbrace u=s e^{i\sigma}$: $-\infty<s<\infty\rbrace$, with $\vert\sigma\vert<\pi/8$:
$$
P(x,y)=\frac{1}{2}\int_{-\infty e^{i\sigma}}^{\infty e^{i\sigma}} e^{-u^4-xu^2+iyu}du.
$$
Define the parameter $\alpha:=\displaystyle{\frac{y}{\sqrt{x}}}$ and $\theta:=\arg x$. After the change of variable $t=u\sqrt{x}=u\vert x\vert^{1/2}e^{i\theta/2}$ we find
\begin{equation}\label{pxypp}
P(x,y)=\displaystyle{\frac{1}{2\sqrt{x}}\int_{-\infty e^{i(\sigma+\theta/2)}}^{\infty e^{i(\sigma+\theta/2)}} e^{-t^2+i\alpha t}e^{-t^4/x^2}dt.}
\end{equation}
From here, an asymptotic expansion of the integral \Eqref{pxypp} for large $\vert x\vert$ can be obtained considering the Taylor expansion of the exponential at the origin:
\begin{equation}\label{expa}
e^{-t^4/x^2}=\sum_{k=0}^{n-1}\frac{(-1)^kt^{4k}}{k!x^{2k}}+r_n\left(\frac{t^4}{x^2}\right),
\end{equation}
where $r_n(t)$ is the Taylor remainder. From the Lagrange's formula for the remainder we have:
$$
 r_n\left(\frac{t^4}{x^2}\right)= r_n\left(s^4e^{i\,4\sigma}\right)=\frac{t^{4n}}{n!x^{2n}}e^{-\xi\,e^{i\,4\sigma}}\,\,\,\,\,\textnormal{with}\,\,\,\,0<\xi<s^4.
$$
Then
\begin{equation}\label{cotar}
\left\vert r_n\left(\frac{t^4}{x^2}\right)\right\vert<\frac{\vert t\vert ^{4n}}{n!\vert x\vert ^{2n}}.
\end{equation}

Replacing expansion \eqref{expa} in the above integral and interchanging sum and integral we find
$$
P(x,y)=\displaystyle{\frac{\sqrt{\pi}\,e^{\frac{-y^2}{4x}}}{2\sqrt{x}}\,\,\left(\sum_{k=0}^{n-1} \frac{(-1)^k}{k!\,x^{2k}}\,Q_{4k}(\alpha)+R_n(x,y)\right)\,\, ,}
$$
with
\begin{equation}\label{erre}
Q_k(\alpha):=\frac{e^{\frac{y^2}{4x}}}{\sqrt{\pi}}\int_{-\infty e^{i(\sigma+\theta/2)}}^{\infty e^{i(\sigma+\theta/2)}} t^{k}e^{-t^2+i\alpha t}dt,\hskip 1cm R_n(x,y):=\frac{e^{\frac{y^2}{4x}}}{\sqrt{\pi}}\int_{-\infty e^{i(\sigma+\theta/2)}}^{\infty e^{i(\sigma+\theta/2)}} r_n\left(\frac{t^4}{x^2}\right)e^{-t^2+i\alpha t}dt.
\end{equation}
These integrals exist when $\vert\arg x+2\sigma\vert<\pi/2$. For any value of $x$ such that $\vert\arg x\vert<3\pi/4$, we can chose $\sigma$ (with $\vert\sigma\vert<\pi/8$) such that these integrals are finite. After some manipulations and using the integral representation of Hermite polynomials [eq.18.10.10, \cite{nisthermite}]

$$\mathop{H_{n}\/}\nolimits\!\left(x\right)=\frac{(-2i)^{n}e^{x^{2}}}{\sqrt{\pi}}\int_{-\infty}^{\infty}e^{-t^{2}}t^{n}e^{2ixt}dt $$

we find

%PONER UNA REFERENCIA CON LA REPRESENTACION INTEGRAL DE LOS HERMITES

$$
Q_{4k}(\alpha)=2^{-4k}\,H_{4k}\left(\frac{\alpha}{2}\right)\,\,\,,\,\,\,\forall k=0,1,2,\ldots
$$

Moreover, introducing \eqref{cotar} in \eqref{erre} we find
$$
\vert R_n(x,y)\vert\le\,\frac{e^{\frac{\vert y\vert^2}{4\vert x\vert}\left(\frac{\cos ^2(\arg y-\arg x-\sigma)}{\cos (2 \sigma+\arg x)}\right)}}{n!\,\vert 4 x\vert^{2n}\cos^{2n+1/2}(2\sigma+\arg x)}\,H_{4n}\left(i\,\frac{\vert y \vert\sin(\sigma+\arg y)}{2\sqrt{\vert x \vert\cos(2\sigma+\arg x)}}\right).
$$

It is clear that $R_n(x,y)=\mathcal{O}(x^{-2n})$ as $\vert x\vert\to\infty$ uniformly in $y$ with bounded $y/\sqrt{x}$ and $\vert\arg x\vert<3\pi/4$.
Therefore, we have obtained that
\begin{equation}\label{expandos}
P(x,y)=\displaystyle{\frac{\sqrt{\pi}\,e^{\frac{-y^2}{4x}}}{2\sqrt{x}}\,\left(\sum_{k=0}^{n-1} \frac{(-1)^k}{k!(4x)^{2k}}H_{4k}\left(\frac{y}{2\sqrt{x}}\right)+R_n(x,y)\right)\,\,,}
\end{equation}
is an asymptotic expansion of $P(x,y)$ when $\vert x\vert\to\infty$ uniformly in $y$ with bounded $y/\sqrt{x}$ and $\vert\arg x\vert<3\pi/4$.

This expansion is just the expansion (3.2) derived in \cite{paris}. The derivation of (3.2) is much more cumbersome as it is obtained from a contour integral representation of the Pearcey integral in terms of a parabolic cylinder function. An error bound for the remainder is not provided there. On the other hand, an asymptotic expansion for $\vert\arg (-x)\vert<\pi/4$ is given in \cite{paris}.

\section{An asymptotic expansion for large $\vert x\vert$ and $\vert y\vert$}

Define the parameter $\gamma:=\displaystyle{\frac{y}{2x}}$. We can write the Pearcey integral \Eqref{pearcy} in the form
\begin{equation}\label{pxy}
P(x,y)=\displaystyle{\frac{1}{2}\int_{-\infty}^{\infty} e^{x f(t)-t^4}dt, }
\end{equation}
with phase function $f(t):=2i\gamma t-t^2$. The unique saddle point of this phase function is the point
$t_0=i\gamma$. From \cite{nico} we know that the asymptotic analysis of this integral does not require the application of the standard saddle point method, but just a change of variable. Using the Cauchy's residue theorem we find that we can replace the integration path $(-\infty,\infty)$ in \Eqref{pxy} by any straight $\Gamma:=\lbrace i\gamma+t$; $t=r e^{i\gamma}$: $-\infty<r<\infty\rbrace$, with $\vert\sigma\vert<\pi/8$:
\begin{equation}\label{pxyp}
P(x,y)=\displaystyle{\frac{e^{-\gamma^2 x}}{2}\int_{-\infty }^{\infty } e^{-xt^2}e^{-(i\gamma+t)^4}dt.}
\end{equation}
We write the above integral in the form
\begin{equation}\label{pepe}
P(x,y)=\displaystyle{\frac{e^{-\alpha^2 x-\alpha^4}}{2}\int_{-\infty e^{i\sigma}}^{\infty e^{i\sigma}} e^{-(x-6\alpha^2)t^2}h(t)dt, }
\end{equation}
with
$$
h(t):=e^{4i\gamma^3 t-4i\gamma t^3-t^4}.
$$
From here, an asymptotic expansion of the integral \Eqref{pxyp} for large $\vert x\vert$ can be derived from Watson's lemma. From the Taylor expansion of the exponential at the origin we have that
$$
h(t)=\sum_{n=0}^\infty A_n(\gamma)t^n, \hskip 1cm A_n(\gamma):=\sum_{k=0}^{\lfloor n/4\rfloor}\sum_{j=0}^{\lfloor(n-4k)/3\rfloor}\frac{(4i\gamma^3)^{n-4k-3j}(-4i\gamma)^j(-1)^k}{k!j!(n-4k-3j)!}.
$$
Replacing this expansion in the integral \eqref{pepe} and interchanging sum and integral we find that
$$
P(x,y)\sim\displaystyle{\frac{e^{-y^2/(4x)-y^4/(16x^4)}}{2}\sum_{n=0}^\infty A_{2n}(\gamma)P_n(\gamma,x),}
$$
with
$$
P_n(\gamma,x):=\int_{-\infty }^{\infty } t^{2n}e^{-(x-6\gamma^2)t^2}dt, \hskip 2cm t=r e^{i\gamma}, \hskip 5mm -\infty<r<\infty.
$$
This integral exists when $\vert\arg(x-6t=r e^{i\gamma}$: $-\infty<r<\infty^2)+2\sigma\vert<\pi/2$. For any value of $x$ and $y$ such that $\vert\arg(x-3y^2/(2x^2))\vert<3\pi/4$, we can chose $\sigma$ (with $\vert\sigma\vert<\pi/8$) such that this integral is finite:
$$
P_n(\gamma,x)=\frac{\Gamma(n+1/2)}{(x-6\gamma^2)^{n+1/2}}.
$$
Therefore, we have obtained that, when $\vert x\vert\to\infty$,
\begin{equation}\label{expantres}
P(x,y)\sim\displaystyle{\frac{e^{-y^2/(4x)-y^4/(16x^4)}}{2}\sqrt{\frac{2x^2}{2x^3-3y^2}}\sum_{n=0}^\infty A_{2n}\left(\frac{y}{2x}\right)\Gamma(n+1/2)\left(\frac{2x^2}{2x^3-3y^2}\right)^n,}
\end{equation}
uniformly in $y$ with bounded $y/x$ and $\vert\arg(x-3y^2/(2x^2))\vert<3\pi/4$.

\vspace{0.5cm}

\section{Final remarks and numerical experiments}

Figure 1 summarizes the regions for $(\vert x\vert,\vert y\vert)$ covered by the above mentioned algorithms (except \cite{pariskaminski} that is more sophisticated and valid only for real $(x,y)$) and by the three algorithms derived below in this paper.

%\begin{figure}
\begin{center}
\includegraphics[width=0.85\textwidth]{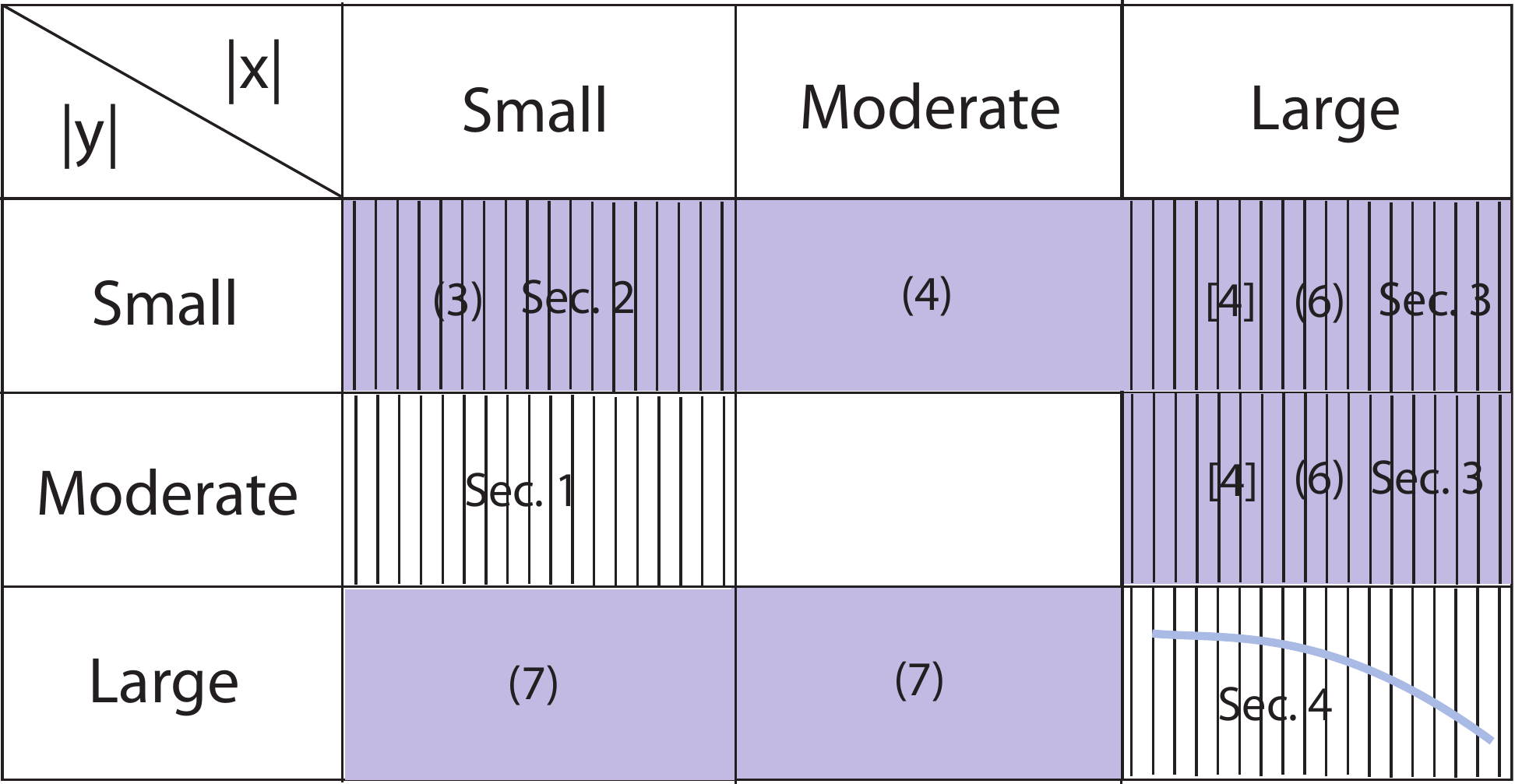}
\end{center}
Figure 1:{\small In blue, the regions of $(\vert x\vert,\vert y\vert)$ for which there are analytic formulas for the computation of $P(x,y)$ available in the literature and detailed in the introduction section. The regions for which the new formulas derived here are useful are dashed. The curve in the large $\vert x\vert$ and $\vert y\vert$ box indicates the caustic.}
%\label{fig:path1}
%\end{figure}

\vspace{0.5cm}

The expansion given in Section 2 is new and valid for any complex values of $x$ and $y$. The expansion given in Section 3 is just the one given in eq. (3.2) of \cite{paris}, although derived in a more straightforward way that permits an accurate error bound for the remainder. It is valid for $\vert\arg x\vert<3\pi/4$ and any complex $y$. The expansion given in Section 4 is also new and permits a wider region for $(x,y)$ than the expansion given in \cite{kaminski}. The expansion given in \cite{kaminski} is restricted to the vicinity of the caustic $8x^3-27y^2=0$, whereas the one derived in Section 4 above is valid for $\vert\arg(x-3y^2/(2x^2)\vert<3\pi/4$.

Finally, we illustrate the accuracy of the above expansions. In the following tables we show the relative error obtained with these algorithms for several values of $(x,y)$ and different orders $n$ of the approximation. As we do not have at our disposal the exact value of the Pearcey integral, we have taken the numerical integration obtained with the program {\it Mathematica} (with {\it working precision}$=200$) as the exact value of $P(x,y)$.

\begin{table}[H]
\begin{center}
{\small{\begin{tabular}{lllllll}
\hline
\hskip 3mm$(x,y)$ & \multicolumn{6}{c}{$n$} \\
 & 0 & 1 & 2 & 3 & 4 & 5\\ \hline
$(1,i)$  &0.392107 & 0.1672 & 0.063132 & 0.021281 & 0.006514 & 0.001836 \\
$(e^{i\pi/4},1)$  &0.257819 & 0.079558 & 0.022810 & 0.005967 & 0.001428 & 0.000314 \\
$(-1/5,i)$  &0.080452 & 0.006563 & 0.000477 & 0.000031 & $1.868\cdot 10^{-6}$ & $1.033\cdot 10^{-7}$ \\
$(1/10-i/8,-2)$  & 0.004970 & 0.002162 & 0.000182 & 0.000011 & $5.876\cdot 10^{-7}$ &$2.719\cdot 10^{-8}$ \\
$(i/20,2i)$ &0.027790 & 0.000675 & 0.000013 & $2.501\cdot 10^-7 $& $4.042\cdot 10^-9$ &$5.976\cdot 10^-11$ \\
\hline
\end{tabular}}}
\end{center}
\caption{Relative error for several values of $(x,y)$ and the number of terms $n$ of the approximation given by \Eqref{expanuno}.}
\label{tab:er1}
\end{table}
%\medpar

\begin{table}[H]
\begin{center}
{\small{\begin{tabular}{lllllll}\hline
\hskip 3mm$(x,y)$ & \multicolumn{6}{c}{$n$} \\
& 0 & 1 & 2 & 3 & 4 & 5\\ \hline
$(5,2i)$ & 0.04689 & 0.011308 & 0.004835 & 0.002943 & 0.002317 & 0.002233 \\
$(10e^{i\pi/4},-1)$  &0.006956 & 0.000281 & 0.000021 & $2.441\cdot 10^{-6}$ & $3.675\cdot 10^{-7}$ &$ 6.934\cdot 10^{-8}$ \\
$(5+10i,1-i)$ &0.007110 & 0.000291& 0.000022& $2.505\cdot 10^{-6}$ & $3.690\cdot 10^{-7}$ & $6.722\cdot 10^{-8}$\\
$(20,1)$ &0.001766 & 0.000018 & $3.514\cdot 10^{-7}$ & $1.002\cdot 10^{-8}$ & $3.772\cdot 10^{-10}$ & $1.846\cdot 10^{-11}$ \\
$(30i,-i)$  &0.000837 & $4.092\cdot 10^{-6}$ & $3.775\cdot 10^{-8}$ & $5.147\cdot 10^{-10}$ & $9.215\cdot 10^{-12}$ & $3.221\cdot 10^{-13}$ \\
$(100i,2-i)$ &0.000078& $3.553\cdot 10^{-8}$ & $3.048\cdot 10^{-11}$ & $3.885\cdot 10^{-14}$ & $1.550\cdot 10^{-16}$ & $1.150\cdot 10^{-16}$ \\
\hline
\end{tabular}}}
\end{center}
\caption{Relative error for several values of $(x,y)$ and the number of terms $n$ of the approximation given by \Eqref{expandos}.}
\label{tab:er2}
\end{table}
%\medpar

%\multicolumn{4}{c}{Non-logarithmic case}\\
%\multicolumn{4}{c}{$a=-0.25$, $b=b'=1$, $c=c'=1.5$, $y=0.5$}

%\begin{figure}
%\begin{center}
%\includegraphics[width=0.35\textwidth]{graf1.eps}\qquad\includegraphics[width=0.35\textwidth]{graf2.eps}
%\end{center}
%\caption{\small Plot of the exact solution $y(x)$ of Example 1 (dashed) and the approximations $y_{1}(x)$ (red), $y_2(x)$ (green) and $y_3(x)$ (blue).}
%\label{fig:ejem3}
%\end{figure}
%\hfill$\square$

\begin{table}[H]
\begin{center}
{\small{\begin{tabular}{lllllll}\hline
\hskip 3mm$(x,y)$ & \multicolumn{6}{c}{$n$} \\
& 0 & 1 & 2 & 3 & 4 & 5\\ \hline
$(10,3i)$ & 0.007738 & 0.007732 & 0.000532 & 0.000353 & 0.000071 & 0.000031 \\
$(20e^{5i\pi/8},10)$  &0.002036 & 0.002008 & 0.000139 & 0.000027 & $5.602\cdot 10^{-6}$ & $7.790\cdot 10^{-7}$ \\
$(50i,20e^{i\pi/4})$ &0.000295 & 0.000185 &$4.224\cdot 10^{-6}$ & $5.134\cdot 10^{-7}$ & $2.801\cdot 10^{-8}$ &$1.784\cdot 10^{-9}$ \\
$(100,20i)$ &0.000074 & 0.000073 & $1.820\cdot 10^{-7}$ & $3.237\cdot 10^{-8}$ & $2.616\cdot 10^{-10}$ & $2.618\cdot 10^{-11}$ \\
$(200i,5)$  &0.000018 &0.000017 &$ 2.062\cdot 10^{-9}$ &$2.051\cdot 10^{-9}$ & $4.274\cdot 10^{-13}$ & $4.228\cdot 10^{-13} $\\
\hline
\end{tabular}}}
\end{center}
\caption{Relative error for several values of $(x,y)$ and the number of terms $n$ of the approximation given by \Eqref{expantres}.}
\label{tab:er3}
\end{table}
%\medpar

%\pagebreak
%
\section{Acknowledgments}

This research was supported by the Spanish \emph{Ministry of "Econom\'{\i}a y Competitividad"}, project MTM2014-52859. The {\it Universidad P\'ublica de Navarra} is acknowledged by its financial support.

\footnotesize{
}
\end{document}